\documentclass{article}
\usepackage{amsfonts}
\usepackage{amsmath}

\newtheorem{theorem}{Theorem}
\newtheorem{acknowledgement}[theorem]{Acknowledgement}

\newtheorem{corollary}[theorem]{Corollary}

\newtheorem{definition}[theorem]{Definition}

\newtheorem{lemma}[theorem]{Lemma}

\newtheorem{proposition}[theorem]{Proposition}
\newtheorem{remark}[theorem]{Remark}

\input{tcilatex}

\begin{document}

\title{Projective relations for m-th root metric spaces}
\author{Nicoleta Brinzei}
\date{}
\maketitle

\begin{abstract}
For Finsler spaces $(M,F)$ endowed with m-th root metrics, we provide
necessary and sufficient conditions in which they are projectively flat, or
projectively related to Berwald/Riemann spaces. We also give a specific
characterization for m-th root metrics spaces of Landsberg and of Berwald
type.
\end{abstract}



\section{Introduction}

Finsler spaces with m-th root metric, defined by the fundamental function%
\begin{equation}
F=\sqrt[m]{a_{i_{1}\dots i_{m}}(x)y^{i_{1}}y^{i_{2}}\dots y^{i_{m}}}
\label{F}
\end{equation}%
have been introduced and studied in 1979 by Shimada, \cite{Shimada}, and, in
the case $m=3$, by Matsumoto and Okubo, \cite{Mat}. Most recently, they are
taken into consideration by physicists as subject for a possible model of
space-time, \cite{1}, \cite{Pavlov2}, \cite{Siparov},\cite{Garasko}, \cite%
{Bogoslovsky}.

For quartic metrics ($m=4$), a study of geodesics and of the related
geometrical objects is made by S. Lebedev, \cite{Leb}, respectively, by V.
Balan, S. Lebedev and the author, \cite{BBL}. Also, Einstein equations for
some relativistic models relying on such metrics are studied by V. Balan and
the author in two papers, \cite{Mangalia}, \cite{Cairo1}. Tensorial
connections for such spaces have been recently studied by L. Tamassy, \cite%
{Tamassy}.

\bigskip

In the following, we shall give answers to the following questions:

1) In what conditions an m-th root metric space has rectilinear geodesics?

2) When do geodesics of such a space coincide with those of a Riemannian one?

3) In what conditions an m-root metric space is of Landsberg or of Berwald
type?

\bigskip

The main idea is to use the Lagrange-type metric tensor%
\begin{equation}
h_{ij}=a_{ij0...0},  \label{flag}
\end{equation}%
(which is nondegenerate), called in the following, the \textit{polynomial
flag metric }(V. Balan, \cite{Mangalia}), whose coefficients $h_{ij}$ are
all polynomial functions in the directional variables $y^{i}$, instead of
metric tensors traditionally used in the study of such spaces, namely, the
usual Finsler metric, or Shimada's $a_{ij}=\dfrac{h_{ij}}{F^{m-2}}$ (which
both contain $m$-th roots in their expressions). The expressions of
corresponding geometrical objects and algorithms are thus simplified.

\bigskip

In Section 2, we express the equations of geodesics and the related
geometrical objects (canonical nonlinear connection, Berwald linear
connection etc.) in terms of the polynomial metric $h_{ij}.$

Section 3 is dedicated to m-th root metric spaces of Landsberg,
respectively, of Berwald type; we show that the characterizing conditions
written in terms of the polynomial flag metric $h_{ij}$ look in similar way
to their correspondents in terms of the usual Finsler metric, namely:

\begin{itemize}
\item that the Berwald connection should be $h$-metrical w.r.t. $h_{ij}$ -
for Landsberg spaces;

\item that $h_{ijk|l}=0$ (with respect to the canonical metrical connection
of $h_{ij}$) - for Berwald spaces.
\end{itemize}

In the next section, we make a short review of (known) results regarding:
projective relations between Finsler spaces, projective flatness, Finsler
spaces projective to Riemannian ones, \cite{Shen}, \cite{Bacso}, \cite%
{Vincze}.

In Sections 5 and 6, we find the conditions upon the coefficients $%
a_{i_{1}\dots i_{m}}(x)$ such that an m-th root metric space should be
projectively flat, respectively, projective to a Riemannian space. By using
the polynomial flag metric $h,$ the determining algorithms reduce to simply
identifying the coefficients of some polynomials.

In the last section, we analyze transformations 
\begin{equation*}
\bar{F}^{m^{\prime }}=\alpha (x,y)F^{m}
\end{equation*}%
and determine the conditions which make them become projective. As a
consequence, we determine classes of Riemann-projective and projectively
flat m-th root metrics.

As a remark, m-root metric functions $F$ in (\ref{F}) are not necessarily
positive definite; still, under the assumption of nondegeneracy, we use the
term of "Finsler" metric, \cite{Shimada}.

\bigskip

\section{m-th root metric spaces}

Let $M^{n}$ be a differentiable manifold of dimension $n$ and class $%
\mathcal{C}^{\infty },$ $TM$ its tangent bundle and $(x^{i},y^{i})$ the
coordinates in a local chart on $TM$. Let $F$ be the following function on $%
M,$ \cite{Shimada}:%
\begin{equation}
F=\sqrt[m]{a_{i_{1}\dots i_{m}}(x)y^{i_{1}}y^{i_{2}}\dots y^{i_{m}}}
\label{1.1}
\end{equation}%
(with $a_{i_{1}\dots i_{m}}$ symmetric in all its indices).

In the following, for a function $f=f(x,y),$ we shall denote by " $,"$ and "~%
$\cdot $~" the partial derivatives w.r.t. $x$ and $y,$ respectively. Also,
if $N$ is a nonlinear connection on $TM,$ we denote by " ; " its associate
covariant derivative.%
\begin{equation*}
f_{;l}=\frac{\delta f}{\delta x^{l}}=\frac{\partial f}{\partial x^{l}}%
-N_{~l}^{r}\frac{\partial f}{\partial y^{r}},\;f\in \mathcal{F}(TM).
\end{equation*}

Let $T$ denote the $m$-th power of $F:$%
\begin{equation}
T=F^{m}=a_{i_{1}\dots i_{m}}(x)y^{i_{1}}y^{i_{2}}\dots y^{i_{m}}.  \label{1}
\end{equation}%
For the $y$-derivatives of $F$ and $T$ we shall omit the dot: 
\begin{equation*}
T_{i}=\frac{\partial T}{\partial y^{i}}=T_{\cdot i},T_{ij}=T_{\cdot
ij},F_{i}=F_{\cdot i},F_{ij}=F_{\cdot ij}
\end{equation*}%
and we denote by null index transvection by $y$ (e.g., $L_{i0}=L_{ij}y^{j})$.

Let%
\begin{equation}
g_{ij}=\frac{1}{2}\frac{\partial ^{2}F^{2}}{\partial y^{i}\partial y^{j}}
\label{usual}
\end{equation}%
the \textit{usual Finsler metric} determined by $F$.

We denote by 
\begin{equation}
h_{ij}=\frac{1}{m(m-1)}T_{ij}=a_{ij0...0},  \label{h}
\end{equation}%
the \textit{polynomial (flag) metric} attached to $F.$

Then we have: 
\begin{equation}
h_{i0}=\frac{1}{m}T_{i}.  \label{2}
\end{equation}%
The normalized supporting element $l_{i}=F_{i}$ can be written in terms of $%
h $ as: 
\begin{equation}
l_{i}=\frac{h_{i0}}{F^{m-1}}.  \label{3}
\end{equation}

\begin{remark}
In \cite{Shimada}, \cite{Mangalia}, one uses the homogenized version%
\begin{equation}
a_{ij}=\dfrac{1}{F^{m-2}}h_{ij}.  \label{a}
\end{equation}
\end{remark}

\bigskip

The link between $h_{ij}$ and the usual Finsler metric $g_{ij}$ is given by:%
\begin{equation}
h_{ij}=F^{m-2}(\dfrac{1}{m-1}g_{ij}+\dfrac{m-2}{m-1}l_{i}l_{j}).  \label{1.2}
\end{equation}%
Taking into account that $g_{ij}$ is nondegenerate on $TM\backslash \{0\}$,
it follows that $h_{ij}$ itself is nondegenerate on $TM\backslash \{0\}$;
its inverse matrix is:%
\begin{equation*}
h^{ij}=\dfrac{1}{F^{m-2}}\{(m-1)g^{ij}-(m-2)l^{i}l^{j}\}.
\end{equation*}

In other words, $h$ is a Lagrange metric, \cite{Miron}, on $M$.

Moreover, if $g$ is positive-definite, then so is the polynomial metric $h.$

The proof of both nondegeneracy and positive-definiteness, as well as the
computation of the inverse matrix $h^{ij}$ rely on \ref{1.2} and the
following lemma (\cite{BBL}):

\begin{lemma}
Consider the following family of (0,2)-tensor fields:%
\begin{equation*}
\Theta _{ij}=\lambda g_{ij}+\mu l_{i}l_{j},~\ \ \lambda ,\mu \in \mathcal{F}%
(M).
\end{equation*}%
Denote by $g^{ij}$ the dual of $g_{ij}.$ Then

a) $\Theta _{ij}$ is non-degenerate for $\lambda (\lambda +\mu )\not=0$ on $%
TM;$

b) The dual of $\Theta _{ij}$ is%
\begin{equation*}
\Theta ^{ij}=\dfrac{1}{\lambda }g^{ij}+\dfrac{-\mu }{\lambda (\lambda +\mu
)F^{2}}y^{i}y^{j};
\end{equation*}

c) the determinant of $\Theta $ is $\det \Theta =\lambda ^{n-1}(\lambda +\mu
)\det g.$
\end{lemma}

By expressing the Euler-Lagrange equations attached to $F$ in terms of $%
h_{ij},$ we get

\begin{proposition}
The equations of unit-speed geodesics $t\mapsto (x^{i}(t))$ of an m-th root
metric space can be expressed in terms of the polynomial metric $h_{ij}$ as%
\begin{equation}
\dfrac{d^{2}x^{i}}{dt^{2}}+\gamma
_{~j_{1}...j_{m}}^{i}y^{j_{1}}...y^{j_{m}}=0,  \label{1.3}
\end{equation}%
where $y^{i}=\dfrac{dx^{i}}{dt}$ and the generalized Christoffel symbols $%
\gamma _{j_{1}...j_{m}}^{i}$ are%
\begin{equation}
\gamma _{~j_{1}...j_{m}}^{i}=\dfrac{h^{ip}}{m(m-1)}\left( \underset{%
(j_{1},...,jm)}{\sum }\dfrac{\partial a_{pj_{2}...j_{m}}}{\partial x^{j_{1}}}%
-\dfrac{\partial a_{j_{1}...j_{m}}}{\partial x^{p}}\right) ,  \label{1.4}
\end{equation}
\end{proposition}

with $\underset{(j_{1},...,jm)}{\sum }$ denoting cyclic sum w.r.t. the
involved indices.

Obviously, $\gamma _{~j_{1}...j_{m}}^{i}$ are totally symmetric w.r.t. the
lower indices. Moreover,%
\begin{equation*}
\gamma _{pj_{1}...j_{m}}=h_{ip}\gamma _{j_{1}...j_{m}}^{i}
\end{equation*}%
depend only on $x.$

We notice that%
\begin{equation*}
\gamma _{p00...0}=\gamma _{pj_{1}...j_{m}}y^{j_{1}}...y^{j_{m}}.
\end{equation*}

are polynomial functions (homogeneous of degree $m$) in $y^{i}.$

\bigskip

According to \cite{Miron}, the canonical spray coefficients of $(M,F)$ are%
\begin{equation}
2G^{i}=\gamma _{~00..0}^{i}.  \label{1.5}
\end{equation}

It is also useful to express the canonical (Kern) nonlinear connection and
the Berwald connection of the given Finsler space in terms of $h_{ij}.$
Namely, the coefficients of the nonlinear connection are

\begin{equation}
N_{~j}^{i}=G_{~\cdot j}^{i}=\frac{1}{2}\frac{\partial }{\partial y^{j}}%
(h^{is}\gamma _{s~00\dots 0})=\frac{1}{2}(h_{~\cdot j}^{is}\gamma _{s~0\dots
0}+mh^{is}\gamma _{s~j0\dots 0}).  \label{1.6}
\end{equation}

In the following, we shall always mean by $N$ the Kern connection (\ref{1.6}%
).

Further, the Berwald connection coefficients $G_{~jk}^{i}=G_{~\cdot jk}^{i}$%
, \cite{Anto}, \cite{Shen}, are: 
\begin{equation}
G_{~jk}^{i}=\frac{1}{2}(h_{~\cdot jk}^{is}\gamma _{s~0\dots 0}+mh_{~\cdot
j}^{is}\gamma _{s~k0\dots 0}+mh_{~\cdot k}^{is}\gamma _{s~j0\dots
0}+m(m-1)h^{is}\gamma _{s~jk0\dots 0}).  \label{1.7}
\end{equation}%
By subsequent derivation, we infer that its $hv$-curvature $%
G_{jkl}^{i}:=G_{~jk\cdot l}^{i}$ is given by%
\begin{equation}
\begin{array}{l}
G_{jkl}^{i}=\frac{1}{2}\{h_{~\cdot jkl}^{is}\gamma _{s~0\dots
0}+m\sum_{(j,k,l)}(h_{~\cdot jk}^{is}\gamma _{s~l0\dots 0})+\medskip \\ 
\quad +m(m-1)\sum_{(j,k,l)}(h_{~\cdot j}^{is}\gamma _{s~kl0\dots
0})+m(m-1)(m-2)h^{is}\gamma _{s~jkl0\dots 0}\}.%
\end{array}
\label{1.8}
\end{equation}

\textbf{Notes:} 1) The equations of geodesics (\ref{1.3}) can also be
expressed as%
\begin{equation*}
\dfrac{d^{2}x^{i}}{dt^{2}}+G^{i}=0\Leftrightarrow \dfrac{d^{2}x^{i}}{dt^{2}}%
+N_{~j}^{i}y^{j}=0\Leftrightarrow \dfrac{d^{2}x^{i}}{dt^{2}}%
+L_{~jk}^{i}y^{j}y^{k}=0.
\end{equation*}

2) The Kern nonlinear connection and the Berwald connection $B\Gamma (N)$
depend only on the function $F,$ \textit{not} on the choice of the metric
tensor ($h_{ij},$ $a_{ij}$ or $g_{ij}$). We cannot say the same of the
canonical metrical connection $C\Gamma (N).$

By \textit{canonical metrical connection }$C\Gamma (N)$ attached to a
generalized Lagrange metric $\alpha _{ij}$, \cite{Miron}, we mean the normal
connection $(L_{~jk}^{i},C_{~jk}^{i})$ given by%
\begin{eqnarray*}
L_{~jk}^{i} &=&\dfrac{1}{2}\alpha ^{ih}(\alpha _{hj;k}+\alpha _{hk;j}-\alpha
_{jk;h}) \\
C_{~jk}^{i} &=&\dfrac{1}{2}\alpha ^{ih}(\alpha _{hj\cdot k}+\alpha _{hk\cdot
j}-\alpha _{jk\cdot h}).
\end{eqnarray*}%
In particular, if $N$ is the Kern nonlinear connection and $\alpha _{ij}$ is
a Finsler metric, then $(N_{~j}^{i},L_{~jk}^{i},C_{~jk}^{i})$ is the Cartan
connection.

3) For both $B\Gamma (N)$ and $C\Gamma (N)$ there holds, \cite{Miron}:%
\begin{equation*}
y_{~|k}^{i}=0,
\end{equation*}%
where $_{|}$ denotes the horizontal covariant derivative.

\bigskip

\section{m-th root metric spaces of Landsberg and Berwald type}

A Finsler space $(M,F)$ is called a \textit{Landsberg space} if the Berwald
connection $B\Gamma (N)$ is $h$-metrical. Equivalently, \cite{Anto}, $(M,F)$
is a Landsberg space iff the horizontal coefficients of $B\Gamma (N)$ and $%
C\Gamma (N)$ coincide:%
\begin{equation*}
G_{~jk}^{i}=L_{~jk}^{i},~\ \forall i,j,k=1,...,n.
\end{equation*}

$(M,F)$ is called a \textit{Berwald space} if the coefficients $G_{~jk}^{i}$
are functions of position $x$ alone. Equivalently: $G_{~jkl}^{i}=0.$

There holds the inclusion, \cite{Anto}:

\begin{center}
Berwald spaces$~\subset ~$Landsberg spaces.
\end{center}

In \cite{Shimada}, there is proven that

1. The horizontal coefficients of the canonical metrical connection $%
^{a}C\Gamma (N)$ attached to the homogenized metric $a_{ij}$ (\ref{a})
coincide with those of the Cartan connection (attached to the usual Finsler
metric) $^{\ast }C\Gamma (N)$ of $(M,F);$

2. an m-th root metric space $(M,F)$ is a Berwald space (resp. Landsberg
space) if and only if, w.r.t. $^{a}C\Gamma (N),$ we have $a_{ijk|h}=0$
(resp. $a_{ijk|0}=0,$ where%
\begin{equation*}
a_{ijk}=\dfrac{a_{ijk00...0}}{F^{m-3}}.
\end{equation*}

\bigskip

Now, by using:%
\begin{eqnarray*}
h_{ij} &=&F^{m-2}a_{ij},~F_{;l}=0, \\
h_{ij\cdot k} &=&(m-2)F^{m-3}a_{ijk},
\end{eqnarray*}%
we can immediately see that there holds

\begin{proposition}
The horizontal coefficients $L_{~jk}^{i}$ of the canonical metrical
connection $^{h}C\Gamma (N)$ attached to the polynomial metric $h$ coincide
with those of the Cartan connection $^{\ast }C\Gamma (N)$ of $(M,F).$
\end{proposition}

\begin{corollary}
An m-th root metric space $(M,F)$ is a Berwald space (resp. Landsberg space)
if and only if, w.r.t. the canonical metrical connection $^{h}C\Gamma (N),$
we have $h_{ij\cdot k|l}=0$ (resp. $h_{ij\cdot k|0}=0$).
\end{corollary}

For Landsberg spaces, we can obtain easier conditions if we use the Berwald
connection instead of the canonical metrical one.

Namely, by using (\ref{1.2}), (\ref{3}) and the equalities $y_{~|k}^{i}=0,$
we get:

\begin{theorem}
An m-th root metric space $(M,F)$ is of Landsberg type if and only if, with
respect to the Berwald connection $B\Gamma (N),$ we have:%
\begin{equation}
h_{ij|k}=0,~\ \ i,j,k=1,...,n.  \label{1.9}
\end{equation}
\end{theorem}

\section{Finsler metrics in projective relation}

\textbf{Definition.} (\cite{Anto}) Let $F^{n}=(M^{n},F(x,y))$ and $\bar{F}%
^{n}=(M^{n},\bar{F}(x,y))$ be two Finsler spaces on the same underlying
manifold $M^{n}.$ If any geodesic of $F^{n}$ coincides with a geodesic of $%
\bar{F}^{n}$ as a set of points and vice versa, then the change $%
F\rightarrow \bar{F}$ is called \textit{projective} and $F^{n}$ is said to
be \textit{projective to} $\bar{F}^{n}.$

In the following, we shall expose several classical results referring to
projective changes:

\begin{theorem}
(Knebelmann), \cite{Anto}, \cite{Bacso}: A Finsler space $F^{n}$ is
projective to another Finsler space $\bar{F}^{n}$ if and only if there
exists a 1-homogeneous scalar field $p(x,y)$ obeying%
\begin{equation}
\bar{G}^{i}=G^{i}+p(x,y)y^{i}  \label{3.1}
\end{equation}
\end{theorem}

A very useful characterization of projectivity is also given by R\`{a}pcs%
\`{a}k's theorem :

\begin{theorem}
(R\`{a}pcs\`{a}k), \cite{Anto}): A Finsler space $F^{n}=(M^{n},F(x,y))$ is
projective to $\bar{F}^{n}=(M^{n},\bar{F}(x,y))$ if and only if $\bar{F}$
satisfies one of the following three equations:

\begin{enumerate}
\item $\bar{F}_{;j}-y^{r}\bar{F}_{;r\cdot j}=0;$

\item $\bar{l}_{j;i}-\bar{l}_{i;j}=0;$

\item $\bar{F}_{;i\cdot j}-\bar{F}_{;j\cdot i}=0,$ where ; means the
covariant derivative in $F^{n}.$
\end{enumerate}
\end{theorem}

\begin{definition}
A Finsler space $F^{n}$ is said to be \textit{projectively flat} if its
geodesics are straight lines.
\end{definition}

There holds

\begin{theorem}
\cite{Anto}, \cite{Bacso}: A Finsler space is projectively flat iff it is
projective to a locally Minkovski space.
\end{theorem}

For Berwald-type projectively flat spaces, we have the following
characterization:

\begin{theorem}
,\cite{Bacso}, A Finsler space $F^{n}$ (of dimension $n$) is a projectively
flat Berwald space if and only if it belongs to one of the following classes:

\begin{enumerate}
\item $n\geq 3:$

a) locally Minkovski spaces;

b) Riemannian spaces of constant curvature;

\item $n=2:$

a) locally Minkovski spaces;

b) Riemannian spaces of constant curvature;

c) spaces $F^{2}$ with $F=\dfrac{\beta ^{2}}{\gamma },$ where $\beta $ and $%
\gamma $ are 1-forms.
\end{enumerate}
\end{theorem}

Obviously, in the case of m-th root metric spaces, we cannot have the
situation 2c), which implies that the only m-th root metric spaces which are
both Berwald and projectively flat are either locally Minkovski, either
Riemannian $(m=2)$ of constant curvature.

\bigskip

Another important category of Finsler spaces are \textit{Douglas spaces.}

\begin{definition}
A Finsler space is said to be of Douglas type, or a Douglas space, if%
\begin{equation}
D^{ij}=G^{i}y^{j}-G^{j}y^{i}  \label{3.2}
\end{equation}%
are homogeneous polynomials in $y^{i}$ of degree three.
\end{definition}

A characterization of Douglas spaces is given by

\begin{theorem}
\cite{Bacso}: A Finsler space is of Douglas type if and only if the Douglas
projective tensor%
\begin{equation}
D_{~jkl}^{i}=G_{~jkl}^{i}-G_{jk\cdot l}y^{i}/(n+1)-(G_{jk}\delta
_{l}^{i}+G_{lj}\delta _{k}^{i}+G_{kl}\delta _{j}^{i})/(n+1),  \label{3.3}
\end{equation}%
(where $G_{jk}=G_{~jki}^{i})$ vanishes identically.
\end{theorem}

The Douglas tensor $D_{~jkl}^{i}$ is invariant under projective
transformations.

Both locally Minkovski and Berwald (including Riemannian) spaces are Douglas
spaces.

Here are several results related to Douglas spaces, \cite{Bacso}:

\begin{theorem}
Any Finsler space projective to a Douglas space is itself a Douglas space.
\end{theorem}

\begin{theorem}
$F^{n}$ is a Berwald space iff it is a Landsberg space and a Douglas space.
\end{theorem}

Finally, the most important to us:

\begin{theorem}
For positive-definite Finsler manifolds, there holds the inclusion:

\textbf{projective to Riemann = projective to Berwald }$\subset $ \textbf{%
Douglas.}
\end{theorem}

The equality in the theorem above is due to Szab\'{o}'s theorem of
metrizability of (positive definite) Berwald manifolds, \cite{Shen}, \cite%
{Vincze}, \cite{Bacso}.

\section{Projectively flat m-th root metric spaces}

In order to study projective flatness for m-th root metric spaces, we use
the condition (1) in R\`{a}pcs\`{a}k's theorem:%
\begin{equation*}
\bar{F}_{;j}-y^{r}\bar{F}_{;r\cdot j}=0.
\end{equation*}%
Let, for the moment, $\bar{F}$ denote an m-th root metric, and $F,$ a
locally Minkovski one.

Since $(M^{n},F)$ is locally Minkovski, in a certain coordinate system we
have $N_{~j}^{i}=0,$ this is, $\bar{F}_{;i}=\bar{F}_{,i}$ and the mentioned
condition can be written by using only usual partial derivatives: 
\begin{equation}
\bar{F}_{,j}-y^{r}\bar{F}_{,r\cdot j}=0.  \label{4.1}
\end{equation}

In terms of $T=\bar{F}^{m},$ we obtain that the last condition is equivalent
to%
\begin{equation}
T(T_{,j}-y^{r}T_{,r\cdot j})=y^{r}(\tfrac{1}{m}-1)T_{j}T_{,r}.  \label{4.2}
\end{equation}

But 
\begin{equation}
T_{,j}-y^{r}T_{,r\cdot j}=-m(m-1)\gamma _{j0...0},T_{j}=mh_{j0},T=h_{00};
\label{deriv}
\end{equation}
consequently, we have:

\begin{theorem}
An m-th root metric space is projectively flat if and only if%
\begin{equation*}
mh_{00}\gamma _{j0...0}=h_{j0}y^{r}h_{00,r}.
\end{equation*}
\end{theorem}

Both terms in the above relation are polynomials in $y^{i}.$

By replacing the terms in the above equality with their expressions in $%
a_{i_{1}...i_{m}}$ and by identifying the corresponding coefficients, it
follows

\begin{theorem}
A m-th root metric space is projectively flat if and only if%
\begin{equation}
m\sum a_{ri_{2}...i_{m}}\gamma _{j~j_{1}...j_{m}}=\sum
a_{ji_{2}...i_{m}}a_{j_{1}...j_{m},~r},  \label{4.3}
\end{equation}
\end{theorem}

$\forall j,r,i_{1},...,i_{m},j_{1},...j_{m}=\overline{1,r}$ and the symbol $%
\sum $ means symmetrization with respect to $i_{1},...i_{m},j_{1},...,j_{m}$.

\section{m-th root metric spaces projectively related to Riemannian spaces}

Throughout this section, we shall need to assume that $F$ is properly
Finslerian, i.e., $(g_{ij})$ in (\ref{usual}) is \textit{positive definite}.
As a remark, the results above hold true also for the case of
(nondegenerate) $g_{ij}$ of arbitrary signature.

As said above, by Szab\'{o}'s theorem, a (positive-definite) Finsler space
is projective to a Riemannian space if and only if it is projective to a
Berwald one.

A necessary condition for $F^{n}$ to be Berwald-projective is that $F^{n}$
should be of Douglas type. Consequently, a necessary condition for a m-th
root metric space to be Berwald-projective is that its Douglas tensor should
vanish:%
\begin{equation*}
D_{~jkl}^{i}=0.
\end{equation*}

\bigskip

A condition which is also sufficient is obtained by starting with condition
1 in Rapcsak's theorem.

Let $\Gamma _{~jk}^{i}(x)$ be the Berwald connection coefficients of some
Berwald space $(M,\widetilde{F})$, $N_{~j}^{i}=\Gamma _{~j0}^{i}$ the
corresponding nonlinear connection and let "$~;~$" denote the associated
covariant derivative.

In terms of $T=F^{n},$ the cited condition writes:%
\begin{equation}
mT(T_{;i}-T_{;r\cdot i}y^{r})=(1-m)T_{;r}T_{i}y^{r}.  \label{5.2}
\end{equation}

(\ref{5.2}) is an equality of homogeneous polynomials of degree $2m$ in $%
y^{i}.$

By a direct computation of the involved derivatives, it follows that (\ref%
{5.2}) is equivalent to%
\begin{equation}
mT(T_{,i}-T_{,r\cdot i}y^{r})+(m-1)T_{i}T_{,r}y^{r}=-\Gamma
_{~00}^{i}\{(1-m)T_{i}T_{j}+mTT_{ij}\}.  \label{5.3}
\end{equation}%
Unfortunately, the matrix with the entries%
\begin{equation*}
A_{ij}=(1-m)T_{i}T_{j}+mTT_{ij}
\end{equation*}%
is degenerate, so we cannot use its inverse in order to separate variables
in the right-hand side. Though, taking into account that in both sides we
have polynomials, we can simply identify their coefficients. What we obtain
is a linear system in $\Gamma _{~jk}^{i}$. Namely, we have

\begin{theorem}
The m-th root metric space $F^{n}$ is Riemann-projective iff there exist the
functions $\Gamma _{~jk}^{i}(x)$ which obey (\ref{5.3}) and, w.r.t.
coordinate changes on $M,$ they obey the rules of transformation of the
coefficients of a linear connection.
\end{theorem}

By replacing the derivatives of $T$ in (\ref{5.3}) and identifying the
coefficients, it follows:

\begin{theorem}
The m-th root metric space $F^{n}$ is Riemann-projective iff there exist the
functions $\Gamma _{~jk}^{i}$ which obey

\begin{eqnarray}
&&a_{j_{1}...j_{m}}\gamma _{i~i_{1}...i_{m}}-\dfrac{1}{m}a_{ij_{2}...j_{m}}%
\cdot a_{i_{1}...i_{m}~,~j_{1}}=  \label{5.4} \\
&=&\Gamma
_{~i_{1}j_{1}}^{j}%
\{-a_{ii_{2}...i_{m}}a_{jj_{2}...j_{m}}+a_{iji_{3}...i_{m}}a_{i_{1}i_{2}j_{3}...j_{m}}\}.
\notag
\end{eqnarray}%
and, w.r.t. coordinate changes on $M,$ they obey the rules of transformation
of the coefficients of a linear connection.
\end{theorem}

\section{The transformation $\bar{T}=\protect\alpha (x,y)T$}

This type of transformation allows to reduce the order of the metric, in the
following sense: if $\bar{T}$ defines an $m$-th root metric as in (\ref{1}),
then, we are trying to find some $\alpha $ such that $(M^{n},\bar{T})$ be
projective to a m-th root metric space whose metric has smaller root order $%
(M^{n},T)$ (in particular, to a space of order two, this is, to a Riemann
space, where $T=a_{ij}(x)y^{i}y^{j}$).

Let $T$ and $\bar{T}$ yield two m-th root metrics of order $m^{\prime }$ and 
$m$ respectively, related by%
\begin{equation*}
\bar{T}=\alpha (x,y)T.
\end{equation*}%
Suppose $m^{\prime }<m.$ Then, $\alpha $ is homogeneous of order $%
p=m-m^{\prime }>0$ w.r.t. $y.$ Let again $N_{~j}^{i}$ denote the
coefficients of the Cartan connection attached to $F=T^{1/m^{\prime }}.$

By means of (\ref{5.2}):%
\begin{equation}
\bar{T}_{;i}-\bar{T}_{;r\cdot i}y^{r}=\frac{1-m}{m\bar{T}}\bar{T}_{;r}\bar{T}%
_{i}y^{r}.  \label{6}
\end{equation}%
and of the equality $T_{;i}=0,$ we obtain that $T$ and $\bar{T}$ are
projectively related iff 
\begin{equation}
\alpha _{;i}-\alpha _{;r\cdot i}y^{r}=\frac{1}{m}\frac{T_{i}}{T}\alpha _{;0}+%
\frac{1-m}{m}\frac{\alpha _{i}}{\alpha }\alpha _{;0}.  \label{7}
\end{equation}

\textbf{Remark:} If $\alpha =\alpha (x)$ (conformal transformations), (\ref%
{7}) is written as:%
\begin{equation*}
\alpha _{,i}=\dfrac{1}{m}\dfrac{T_{i}}{T}\alpha _{,r}y^{r},
\end{equation*}%
or%
\begin{equation*}
\alpha _{,r}(\delta _{i}^{r}-\dfrac{1}{m}\dfrac{T_{i}}{T}y^{r})=0.
\end{equation*}

\bigskip

We now return to $\alpha =\alpha (x,y).$ Let us notice that any $\alpha $
with vanishing covariant derivative $\alpha _{;~i}=0$ obeys (\ref{7}). As a
consequence, we get

\begin{proposition}
Let $\gamma =(\gamma _{ij}(x))$ be a Riemannian metric on $M$ and $\alpha
=\alpha (x,y)$ a polynomial function in $y$ (of degree $m>0$) with the
property%
\begin{equation}
\alpha _{;i}=0,  \label{prop}
\end{equation}%
where the covariant derivative is taken w.r.t. the canonical linear
connection (\ref{1.6}) arising from $\gamma $. Then, the m-th root metric
space $(M,F)$ with%
\begin{equation}
F^{m+2}=\bar{T}=\alpha \gamma  \label{transf}
\end{equation}%
is projective to the Riemannian space $(M,\gamma ).$
\end{proposition}

In particular, above we consider a projectively flat Riemannian metric $%
\gamma ,$ we get

\begin{proposition}
If $\gamma $ in (\ref{transf}) is a projectively flat Riemannian metric,
then $\bar{T}$ in (\ref{transf}), \ref{prop}) is a projectively flat m-th
root metric.
\end{proposition}

\begin{acknowledgement}
The work was supported by Romanian Academy grant no. 5/5.02.2008.
\end{acknowledgement}

\bigskip

\bigskip

Author's adress:\newline
\begin{tabular}[t]{l}
Nicoleta Br\^{\i}nzei \\ 
"Transilvania" University, \\ 
50 Iuliu Maniu Street, \\ 
RO-500091, Brasov \\ 
e-mail: nico.brinzei@rdslink.ro%
\end{tabular}

\end{document}